\documentclass[12pt,a4paper]{article}
\usepackage[utf8]{inputenc}

\usepackage{amssymb}
\usepackage{amsmath}
\usepackage{amsthm}

\usepackage{mathptmx}       
\usepackage{helvet}         
\usepackage{courier}        
\usepackage{type1cm}        
\usepackage{makeidx}         
\usepackage{graphicx}        
\usepackage[bottom]{footmisc}

\usepackage{amsfonts}
\usepackage{mathtools}

\theoremstyle{plain}
\newtheorem{theorem}{Theorem}
\newtheorem{proposition}[theorem]{Proposition}

\theoremstyle{definition}
\newtheorem{definition}[theorem]{Definition}
\theoremstyle{remark}

\makeindex             

\DeclareMathOperator{\supp}{supp}
\DeclareMathOperator{\id}{id}
\DeclareMathOperator{\Hom}{Hom}

\newcommand{\cU}{\mathcal{U}}

\newcommand{\cH}{\mathcal{H}}
\newcommand{\cB}{\mathcal{B}}

\newcommand{\cD}{\mathcal{D}}
\newcommand{\cA}{\mathcal{A}}
\newcommand{\cE}{\mathcal{E}}
\newcommand{\cG}{\mathcal{G}}
\newcommand{\cL}{\mathcal{L}}

\newcommand{\bR}{\mathbb{R}}
\newcommand{\e}{\varepsilon}
\newcommand{\coleq}{\coloneqq}

\newcommand{\cM}{\mathcal{M}}
\newcommand{\cN}{\mathcal{N}}
\newcommand{\abso}[1]{\lvert#1\rvert}
\newcommand{\norm}[1]{\lVert#1\rVert}
\newcommand{\bN}{\mathbb{N}}
\newcommand{\ud}{\mathrm{d}}
\newcommand{\csub}{\subset\subset}
\newcommand{\csup}{\supset\supset}
\newcommand{\cP}{\mathcal{P}}
\newcommand{\cI}{\mathcal{I}}

\newcommand{\unterstrich}{{-}}
\newcommand{\poly}{\cP}
\newcommand{\poli}{\cI}
\newcommand{\pd}{\partial}

\begin{document}

\title{On association in Colombeau algebras without asymptotics}
\author{Eduard A.\ Nigsch\footnote{University of Vienna, Faculty of Mathematics, Oskar-Morgenstern-Platz 1, 1090 Vienna, Austria, eduard.nigsch@univie.ac.at}}
%
%
\maketitle

\abstract{A recent variant of Colombeau algebras does not employ asymptotic estimates for its definition. We discuss how the concept of association with distributions transfers to this setting and why it still needs to be based on asymptotics.}

\section{Introduction}\label{sec:1}

Colombeau algebras of nonlinear generalized functions are given by factor spaces of suitable spaces of moderate and negligible functions (we refer to the standard literature \cite{ColNew,ColElem,GKOS,MOBook,Biagioni,zbMATH01226424}). Usually, whether a function $R$ in a certain basic space is moderate or negligible is decided by evaluating it on a suitable test object $(\vec\varphi_\e)_{\e \in (0,1]}$ (essentially an approximate identity) and examining the resulting asymtotic behaviour as $\e \to 0$, see \cite{specfull} for a unifying discussion. In special Colombeau algebras, moreover, evaluation on $\vec\varphi_\e$ is already built into the basic space in the sense that the embedded image of a distribution $u \in \cD'$ is given by the net $\langle u, \vec\varphi_\e(x) \rangle$ of smooth functions in $x$. Hence, the basic space for the special algebra on an open subset $\Omega \subseteq \bR^n$ is given by all nets $(u_\e)_\e \in C^\infty(\Omega)^{(0,1]}$, and the growth or vanishing rate of derivatives of $u_\e(x)$ on compact sets determines whether this representative is moderate or negligible, respectively.

Two observations are in order: first, the properties of a given Colombeau algebra directly depend on the test objects $(\vec\varphi_\e)_\e$ used for its definition. Such properties are, for example, diffeomorphism invariance (cf.~\cite{global,papernew}), the possibility to restrict to subspaces (cf.~\cite[Ch.~III, \S 11, p.~100]{MOBook}) but also association properties, as we will see below. Second, so far the choice of test objects had to be made in advance because the basic definitions depend on it.

Recently, a novel construction of Colombeau algebras was given that is closer in spirit to the definition of distributions as those linear functionals on test functions satisfying appropriate continuous seminorm estimates \cite{efree}. Because this change of perspective is similar to that encountered when switching from the sequential approach to distribution theory \cite{zbMATH03421400} to the classical approach based on the theory of locally convex spaces \cite{TD}, these algebras were termed \emph{Colombeau algebras without asymptotics}.

This formulation has several pleasing features. Most importantly, it decouples the definition of the Colombeau algebra itself from the choice of test objects. For this reason, the resulting space is close to being universal in the sense that one has canonical mappings into most of the classically used Colombeau algebras.

However, test objects necessarily reappear in the study of association. While there is no inherent notion of association anymore, we will see that the separation of association tests from the definition of Colombeau algebras makes them more flexible and to a certain degree even arbitrary.

After recalling the definition of a Colombeau algebra without asymptotics in Section \ref{sec:2} and previous notions of association in Section \ref{sec:3} we will discuss association in the new setting in Section \ref{sec:4}.

\section{The Colombeau algebra}\label{sec:2}

In short, the approach of \cite{efree} is based on two steps: first, one extends the domain of distributions from $\cD(\Omega)$ to $C^\infty(\Omega, \cD(\Omega))$ in a natural way, essentially by mapping $u \in \cD'(\Omega)$ to $\id_{C^\infty(\Omega)} \otimes u$. One then still has seminorm estimates of the form $p(u(\vec\varphi)) \le q(\vec\varphi)$, where $p$ and $q$ are continuous seminorms of $C^\infty(\Omega)$ and $C^\infty(\Omega, \cD(\Omega))$, respectively. This step amounts to representing distributions by their regularizations. In a second step one lets go of linearity and replaces the linear estimates in $q(\vec\varphi)$ by polynomial ones. This allows in particular for products to be formed and gives the space of moderate functions. Similarly, a notion of negligible function is obtained by noticing that the prototypical negligible functions $\iota(f)\iota(g)-\iota(fg)$ or $\iota(f) - \sigma(f)$ with $f,g \in C^\infty(\Omega)$ and $\iota,\sigma$ the canonical embeddings, have a uniformly continuous extension to $C^\infty(\Omega, \cE'(\Omega))$ and take the value $0$ if evaluated at the function $\vec\delta \colon x \mapsto \delta(.-x)$ \cite[Lemma 3.1, p.~188]{bigone}.

In the following definition, the sheaves $C^\infty(\unterstrich, \cD(\Omega))$ and $C^\infty(\unterstrich)$ are considered with values in the category of locally convex spaces with smooth mappings in the sense of convenient calculus \cite{KM} as morphisms. The \emph{basic space} of nonlinear generalized functions on $\Omega$ then is the set of sheaf homomorphisms
\[ \cB(\Omega) \coleq \Hom ( C^\infty(\unterstrich, \cD(\Omega)), C^\infty(\unterstrich)). \] 
The embeddings $\iota \colon \cD'(\Omega) \to \cB(\Omega)$ and $\sigma \colon C^\infty(\Omega) \to \cB(\Omega)$ are given by
\begin{alignat*}{2}
 (\iota u)(\vec\varphi)(x) & \coleq \langle u, \vec\varphi(x) \rangle & \qquad &(u \in \cD'(\Omega)) \\
 (\sigma f)(\vec\varphi)(x) & \coleq f(x) & \qquad &(f \in C^\infty(\Omega))
\end{alignat*}
for $\vec\varphi \in C^\infty(U, \cD(\Omega))$ with $U \subseteq \Omega$ open and $x \in U$. For the present discussion we omit the discussion of derivatives \cite[Def.~4]{efree}.

The algebraic structure of moderate and negligible functions is based on the following semirings of polynomials with nonnegative coefficients, $k \in \bN_0$:
\begin{align*}
 \poly_k & \coleq \bR^+ [y_0, \dotsc, y_k], \\
 \poli_k & \coleq \{ \lambda \in \bR^+ [y_0, \dotsc, y_k, z_0, \dotsc, z_k]\ |\ \lambda(y_0, \dotsc, y_k, 0, \dotsc, 0) = 0 \}.
\end{align*}
For $K,L \csub \Omega$, $m,l \in \bN_0$ and $B \subseteq C^\infty(\Omega)$ bounded we set
\begin{alignat*}{2}
\norm{f}_{K,m} & \coleq \sup_{x \in K, \abso{\alpha} \le m} \abso{\pd^\alpha f(x)} & \qquad & (f \in C^\infty(\Omega)), \\
\norm{\vec\varphi}_{K,m; L, l} & \coleq \sup_{\substack{x \in K, \abso{\alpha} \le m \\ y \in L, \abso{\beta} \le l}} \abso{ \pd_x^\alpha \pd_y^\beta \vec\varphi(x)(y) } & \qquad & (\vec\varphi \in C^\infty(\Omega, \cD(\Omega))), \\
\norm{\vec\varphi}_{K,m; B} & \coleq \sup_{\substack{x \in K, \abso{\alpha} \le m \\ f \in B}} \abso{ \langle f(y), \pd_x^\alpha \vec\varphi(x)(y) \rangle } & \qquad & (\vec\varphi \in C^\infty(\Omega, \cE'(\Omega))).
\end{alignat*}
In the following definition, $\cU_x$ denotes the filter base of open neighborhoods of $x$ in $\Omega$.
\begin{definition}\label{def_modnegquot}
An element $R \in \cB(\Omega)$ is called \emph{moderate} if
\begin{gather*}
(\forall x \in \Omega)\ (\exists U \in \cU_x(\Omega))\ (\forall K,L \csub U)\ (\forall m,k \in \bN_0)\\
(\exists c,l \in \bN_0)\ (\exists \lambda \in \poly_k) \ (\forall \vec\varphi_0,\dotsc,\vec\varphi_k \in C^\infty(U, \cD_L(U))):\\
\norm{ \ud^k R(\vec\varphi_0)(\vec\varphi_1,\dotsc,\vec\varphi_k)}_{K, m} \le \lambda ( \norm{\vec\varphi_0}_{K,c; L, l}, \dotsc, \norm{\vec\varphi_k}_{K,c; L, l}).
\end{gather*}
The subset of all moderate elements of $\cB(\Omega)$ is denoted by $\cM(\Omega)$.

An element $R \in \cB(\Omega)$ is called \emph{negligible} if
\begin{gather*}
(\forall x \in \Omega)\ (\exists U \in \cU_x(\Omega))\ (\forall K,L \csub U)\ (\forall m,k \in \bN_0)\ (\exists c,l \in \bN_0)\\
(\exists \lambda \in \poli_k)\ (\exists B \subseteq C^\infty(\Omega)\ \textrm{bounded})\ (\forall \vec\varphi_0, \dotsc, \vec\varphi_k \in C^\infty(U, \cD_L(U))):\\
\norm{\ud^kR(\vec\varphi_0)(\vec\varphi_1,\dotsc,\vec\varphi_k)}_{K, m} \\
\le \lambda ( \norm{\vec\varphi_0}_{K,c; L, l}, \dotsc, \norm{\vec\varphi_k}_{K,c; L, l}, \norm{ \vec\varphi_0 - \vec\delta}_{K, c; B}, \norm{ \vec\varphi_1}_{K, c; B}, \dotsc, \norm{ \vec\varphi_k}_{K, c; B}).
\end{gather*}
The subset of all negligible elements of $\cB(\Omega)$ is denoted by $\cN(\Omega)$.

We set $\cG(\Omega) \coleq \cM(\Omega) / \cN(\Omega)$.
\end{definition}

To connect this definition to those of classical Colombeau algebras think of setting $k=0$ and $\vec\varphi_0(x)(y) = \e^{-n} \varphi(\frac{y-x}{\e})$ for a mollifier $\varphi \in \cD(\bR^n)$. Further details are given in \cite{efree}.

\section{Association in previous contexts}\label{sec:3}

There are canonical mappings from $\cG(\Omega)$ into the algebras $\cG^s(\Omega)$ (the special algebra, cf.~\cite[Section 1.2]{GKOS}), $\cG^e(\Omega)$ (the elementary full algebra, cf.~\cite[Section 1.4]{GKOS}) and $\cG^d(\Omega)$ (the diffeomorphism invariant algebra, cf.~\cite{papernew,bigone}). Because in each of these there is an intrinsic notion of association we first characterize when the image of an element of $\cG(\Omega)$ is associated with a distribution there.

For the special algebra, let the canonical mapping $\Theta^s \colon \cG(\Omega) \to \cG^s(\Omega)$ be given by $(\Theta^s R)_\e(x) \coleq R(\vec\psi_\e)(x)$, where the mollifier $\vec\psi_\e$ used for the embedding into $\cG^s(\Omega)$ is as in \cite{zbMATH06172905}, i.e., $\varphi_\e(x) = \chi ( x \abso{\ln \e})\e^{-n} \rho(x/\e)$ with a cut-off function $\chi$ and a mollifier $\rho$ and $\vec\psi_\e = \varphi_\e(x-y)$.

For the elementary full algebra, the canonical mapping $\Theta^e \colon \cG(\Omega) \to \cG^e(\Omega)$ is given by $(\Theta^e R)(\varphi,x) = R(\vec\varphi)(x)$ for $\varphi \in \cD(\Omega)$ and $x \in \Omega$ with $x + \supp \varphi \subseteq \Omega$, where $\vec\varphi$ is any element of $C^\infty(\Omega, \cD(\Omega))$ such that $\vec\varphi(x)(y)=\varphi(y-x)$ for $y$ in a neighborhood of $x$. Let $\cA_q(\bR^n)$ denote the space of test functions having integral one and vanishing moments of order up to $q$.

Finally, for the diffeomorphism invariant algebra we have a canonical mapping $\Theta^d \colon \cG(\Omega) \to \cG^d(\Omega)$ given by $\Theta^d(R)(\varphi)(x) \coleq R(\vec\varphi)(x)$ with $\vec\varphi(x')\coleq \varphi$ for all $x' \in \Omega$ and $\varphi \in \cD$. Moreover, $S(\Omega)$ denotes the space of test objects for $\cG^d(\Omega)$ \cite[p.~189]{bigone}.

From the respective definitions of association we immediately obtain:

\begin{proposition}\label{prop1}
Let $R \in \cG(\Omega)$ and $u \in \cD'(\Omega)$. Then
\begin{align*}
\Theta^s(R) &\approx u \Longleftrightarrow \forall \psi \in \cD(\Omega): \lim_{\e \to 0} \langle R ( \vec\psi_\e ), \psi \rangle = \langle u, \psi \rangle, \\
\Theta^e(R) &\approx u \Longleftrightarrow \forall \psi \in \cD(\Omega)\ \exists q>0\ \forall \varphi \in \cA_q(\bR^n): \\
&\hphantom{\approx u \Longleftrightarrow} \quad \lim_{\e \to 0} \langle R(S_\e\varphi, .), \psi \rangle = \langle u, \psi \rangle,  \\
\Theta^d(R) &\approx u \Longleftrightarrow \forall \psi \in \cD(\Omega)\ \forall (\vec\varphi_\e)_\e \in S(\Omega): \lim_{\e \to 0} \langle R(\vec\varphi_\e), \psi \rangle  = \langle u, \psi \rangle.
\end{align*}
\end{proposition}

Here, $S_\e(\varphi)(x) \coleq \e^{-n} \varphi(x/\e)$. We list some possible generalizations occurring in the literature:
\begin{itemize}
 \item Strong association requires convergence of order $\e^\beta$ for some $\beta>0$ uniformly for all $\psi$ having support in a given compact set \cite[Def.~1.38, p.~45]{zbMATH01226424}.
 \item $C^k$-association requires convergence to take place in $C^k$ \cite[Def.~3.2.11, p.~287]{GKOS}.
 \item $s$-association (for $s>0$) requires convergence of order $o(\e^s)$ for all $\psi$, while for $D$-$s$-association one in addition takes $\psi$ only from a test function space $D$ \cite[Def.~2.1, p.~92]{zbMATH01226424}.
\item Strong average association replaces the limit in strong association by an averaged limit \cite{zbMATH06749302}.
\end{itemize}

Naturally, all of these can be formulated in the spirit of Proposition \ref{prop1} as well.

\section{Association in the asymptotic-free algebra}\label{sec:4}

A sensible notion of association of $R \in \cB(\Omega)$ with $u \in \cD'(\Omega)$, written $R \approx u$, requires at least the following properties:
\begin{gather}
\label{eq1} \forall u \in \cD'(\Omega): \iota(u) \approx u, \\
\label{eq2} \forall R \in \cN(\Omega): R \approx 0.
\end{gather}
Condition \eqref{eq1}, which ensures minimal compatibility with the distributional world, can be realized by calling $R$ associated to $u$ if $R(\vec\varphi) \to u$ in $\cD'(\Omega)$ as $\vec\varphi$ converges to $\vec \delta: x \mapsto \delta_x$ (i.e., the identity in $\cL(\cD'(\Omega), \cD'(\Omega))$ if we identify kernels and their operators here). Similarly, one can have stronger convergence (e.g.~as in $C^k$-association) of $\iota(u)(\vec\varphi)$ for $u \in \cH$ if one supposes that $\vec\varphi$ converges to the identity on spaces of distributions $\cH \subseteq \cD'(\Omega)$. However, although convergence like $\vec\varphi_\e \to \id$ in $\cL(\cH, \cH)$ gives compatibility with the linear theory, more structure on the test objects seems to be needed to incorporate certain nonlinear effects related to association.
As an example, consider the property
\begin{equation}
x^k \delta^k \approx 0\quad\textrm{ in }\cG^s(\bR),\quad k \in \bN\label{eq5}
\end{equation}
which holds due to
\begin{gather}
 \langle x^k \varphi_\e(x)^k, \psi(x) \rangle = \int x^k \psi(x) \left( \chi(x \abso{\ln \e}) \frac{1}{\e} \rho\left( \frac{x}{\e} \right) \right) ^k\ud x \nonumber \\
 = \e \int z^k \psi(\e z) \chi(\e z \abso{\ln \e})^k \rho(z)\,\ud z = O(\e). \label{eq4}
\end{gather}
A similar calculation holds in the diffeomorphism invariant algebra $\cG^d(\bR)$ of \cite{found}, where the convolution kernel $\vec\psi_\e(x)(y) \coleq \varphi_\e(x-y)$ which is used for the embedding into $\cG^s$ is replaced by $\vec\varphi_\e(x)(y) \coleq \frac{1}{\e} \phi(\e,x)\left(\frac{y}{\e}\right)$ with $\phi \in C^\infty_b(I \times \bR, \cA_q(\bR))$ (cf.~\cite[Definition 7.20]{found}).

Relation \eqref{eq5} is crucially based on the fact that these kernels are obtained by scaling given test functions. Moreover, properties like \eqref{eq4} are essential in applications for calculating associated distributions, as is seen for example in \cite{JVClarke}, where it is shown that the curvature of a conical metric is proportional to the delta distribution at the apex of the cone.

Therefore, we are led to raise the following questions:
\begin{enumerate}
 \item[Q1.] Which (nonlinear) association properties similar to \eqref{eq5} can one expect in general using the convolution kernels of $\cG^s$ or $\cG^d$?
 \item[Q2.] Can the respective association tests be formulated in terms not involving asymptotics?
\end{enumerate}

For condition \eqref{eq2}, which is needed for association to be independent of representatives, we would like to use negligibility to show that $R(\vec\varphi) \to 0$ in $C^m(\Omega)$ ($m=0$ is sufficient) and hence in $\cD'(\Omega)$ if $\vec\varphi \to \vec\delta$ suitably. Negligibility as in Definition \ref{def_modnegquot} implies that for given $K \csub \Omega$, $m \in \bN_0$ and $L \csup K$ we have
\begin{equation}\label{eq3}
\norm{R(\vec\varphi_\e)}_{K,m} \le \lambda ( \norm{\vec\varphi_\e}_{K,c;L,l}, \norm{\vec\varphi_\e-\vec\delta}_{K,c;B})
\end{equation}
whenever $\vec\varphi_\e(x) \in \cD_L(\Omega)$ for $x \in K$. Suppose for simplicity that $\lambda(x,y) = C x^a y^b$ for all $x,y \in \bR_+$, some $C>0$ and  $a,b \in \bN$.
In general, \eqref{eq3} does not necessarily converge to zero even if $\vec\varphi_\e \to \vec\delta$; take for example in dimension $n=1$ a model delta net $\vec\varphi_\e(x)(y) \coleq \e^{-1} \varphi((y-x)/\e)$ where $\varphi \in \cD(\Omega)$ has integral one and, say, vanishing moments of order up to $q$ and nonvanishing $(q+1)$th moment. Then
\[ \norm{\vec\varphi_\e}_{K,c; L, l} = O( \e^{-c-l-1} ) \]
for some constant $c$; however, taking $B = \{ f \}$ with $f(x) \coleq x^{q+1}$ we only obtain
\[ \norm{ \vec\varphi_\e - \vec\delta }_{K,0; B} = \e^{q+1} \abso{ \int z^{q+1} \varphi(z)\,\ud z }. \]
Hence, if $a(c+l+1)-b(q+1)>0$ we cannot conclude that $R(\vec\varphi_\e) \to 0$. We can only do so if sufficiently many moments of $\varphi$ vanish, or more generally, if $\vec\varphi_\e \to \delta$ fast enough. As above, we ask:

\begin{enumerate}
 \item[Q3.] Can the conditions which ensure that \eqref{eq2} holds be formulated in terms not involving asymptotics?
\end{enumerate}

\section{Conclusion}\label{sec:5}

We have seen that if one wants to formulate useful association tests in the asymptotic-free Colombeau algebra one still has to resort to the classically used association tests of full and special Colombeau algebras, which do employ asymptotics. On the one hand this is needed to prove that association does not depend on representatives; on the other hand, this helps in ensuring that association does not only give compatibility with the linear theory, but also is able to handle nonlinear effects as illustrated by the property $x^k \delta^k \approx 0$ ($k \in \bN$), for example. It will be object of further research to investigate whether these association tests necessarily require a formulation in terms of convolution with scaled mollifiers, or whether a more general formulation is possible.

{\bfseries Acknowledgement. }This research was supported by grant P26859 of the Austrian Science Fund (FWF). The author thanks the anonymous referee for comments.

\end{document}